\documentclass[11pt,hyp]{article}
\usepackage{natbib,hyperref}
\usepackage{doi}

\usepackage[leqno]{amsmath}
\usepackage{amssymb,amsthm,graphicx,enumitem,upref}

\usepackage[titletoc, toc, title]{appendix}
\usepackage{breakcites}
\numberwithin{equation}{section}
\theoremstyle{plain}
\newtheorem{theorem}{Theorem}[section]

\theoremstyle{definition} 
\newtheorem{example}[theorem]{Example}
\newtheorem{definition}[theorem]{Definition}
\newtheorem{proposition}[theorem]{Proposition}
\newtheorem{lemma}[theorem]{Lemma}
\newtheorem{corollary}[theorem]{Corollary}

\begin{document}

\title{
A Note on Over- and Under-Representation Among \\
Populations with Normally-Distributed Traits
}

\author{Ronald F. Fox and Theodore P. Hill}

\date{\vspace{-5ex}}  

\maketitle

\begin{abstract}
In every finite mixture of different normal distributions, there will always be exactly one of those distributions that not only is over-represented in the right tail of the mixture, but even completely overwhelms all other subpopulations in the rightmost tails. This property, although not unique to normal distributions, is not shared by other common continuous centrally-symmetric unimodal distributions such as Laplace, nor even by other bell-shaped distributions such as Cauchy (Lorentz) distributions.\end{abstract}

\section{Introduction}

The issue of over- and under-representation of various factions of our society in various roles has become an important subject of study  (e.g., \cite{novob18,  vob19, novob34, novob35, novob36}). A common assumption by experts is that many of the crucial traits in these roles have approximately normal (Gaussian) distributions (e.g., \cite{vob198, vob59, vob191}). The goal of this short note is to present several basic properties of normal distributions that are particularly pertinent to questions of over- and under-representation in the tails. These properties, although not unique to normal distributions, are not shared by other common continuous centrally-symmetric unimodal distributions, even by other bell-shaped distributions such as classical Cauchy (Lorentz) distributions. In particular, in every finite mixture of different normal distributions, there will always be exactly  one of those distributions that not only is over-represented in the right tail but even completely overwhelms all other subpopulations in the rightmost tails. Similar conclusions follow about the left tails, of course, but the emphasis here on the right tails is simply because in practice it is the right tails that receive the most attention - the best athletes and students, for example, not the worst, which are rarely identified.

\section{Normal Distributions in Population Studies}

When population research studies report only the means and standard deviations of their results, the default scientific understanding is that the data are approximately normally distributed, i.e., the distributions in question are close to normal (Gaussian) distributions. 
\begin{figure}[!ht] 
  \centering
  \includegraphics[width=1.0\textwidth]{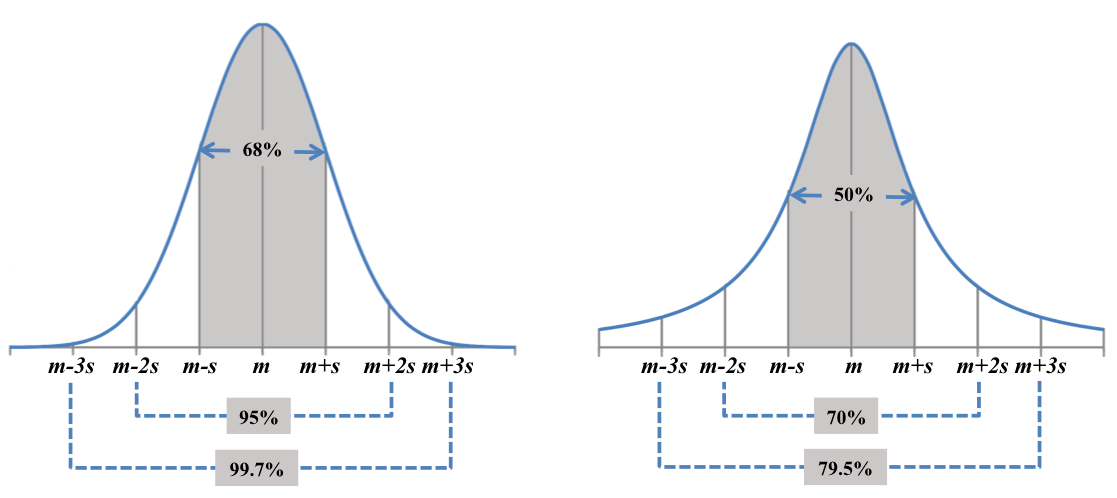}
  \caption{Universal empirical rules for all normal (left) and Cauchy (right) distributions. Here $m$ represents the median, and $s$ represents the distance (called standard deviation for normal distributions, scale parameter for Cauchy) from $m$ to the inflection point. Note that the tails of the normal distribution drop off much more rapidly than those of the Cauchy; for example, the proportion of a normal distribution that is more than 5s (5 standard deviations) out in the right tail is about 6 in one million, whereas for a Cauchy it is about 6 in one hundred. }
  \label{fig:first}
  \end{figure} 
For example, if a research study reports that their data has average value 2 and standard deviation 1, then the usual understanding is that the underlying dataset looks like the diagram in Figure~\ref{fig:first}(left) with $m=2$ and $s=1$, not like the somewhat similar distribution in Figure~\ref{fig:first}(right). The underlying theoretical basis for the assumption of normality in most cases is the remarkable Central Limit Theorem, which says that if the numerical results of independent repetitions of \textit{any experiment} are added, the empirical distribution always approaches a normal distribution.

Two concrete examples among human populations are height and test scores such as those in the College Board Scholastic Aptitude Test (SAT).  There are enormous amounts of data on human height, which is essentially continuous and is very close to being normally distributed \cite[p.~24]{vob191}.  Scores on the SAT, on the other hand, are originally discrete but the distribution ``obtained from a continuized, smoothed frequency distribution of original SAT scores" is a linear transformation of a normal distribution \cite[p.~59]{vob198}.  Thus, since all linear transformations of normal distributions are normal, for all practical purposes the resulting smoothed SAT scores have normal distributions.

The appropriateness of assuming that given data has a normal distribution is often tested using the well-known empirical observation called the ``68\%-95\%-99.7\% rule" of normality. As is illustrated in Figure~\ref{fig:first}(left), in every normal distribution, no matter what the mean and standard deviation (positive square root of variance) are, about 68\% of the values are within one standard deviation of the mean, about 95\% within two standard deviations, and about 99.7\% within three standard deviations (the exact values are irrational numbers, but are easily approximated to any desirable degree of accuracy). The one key property of a continuous centrally-symmetric unimodal distribution that makes it normal is the unique (after rescaling) rate of decrease of its density function away from its mean.  The normal density function, discovered by Gauss in 1809 in connection with his studies of astronomical observation errors, decreases from its mean at a rate  exactly  proportional to $e^{-x^2}$, not to $e^{-x}$ or $x^{-2}$, for example, as is the case for the Laplace and Cauchy distributions, respectively.

The Cauchy distribution, for instance, which sees widespread application in physics, also has a continuous centrally-symmetric unimodal bell-shaped density like the normal distribution, but the Cauchy distribution has an undefined mean and variance and has a different rate of decrease from the median that satisfies a different empirical rule. Figure~\ref{fig:first}(right) illustrates the corresponding 50\%-70\%-79.5\% empirical rule for Cauchy distributions - exactly 50\% of the values are within one scale parameter (distance from median to inflection point of the density curve), and about 70\% within two and 79.5\% within three scale parameters. Thus it is very easy in practice to distinguish between these two similar-looking common bell-shaped distributions.

Clearly the density functions of every two different normal distributions intersect either in exactly one or in exactly two distinct points. Thus the density function of one of those two distributions is strictly larger than that of the other at all points greater than the larger of the two intersection points (or the unique one, if there is only one). This, in turn, implies that the proportion of that distribution from that point on is strictly larger than that of the other distribution from that point on, and thus that this distribution will be over-represented in the right tail. In fact, much stronger conclusions hold, and it is the goal of this note to present a few of these. 

First, however, as will be seen in the next proposition, these points of intersection are easily obtained explicitly, and are recorded here for convenience; the proof follows from basic algebra.

\begin{proposition} \label{prop4feb-1}
The density functions of two different normal distributions with means $m_1$ and $m_2$ and positive standard deviations $\sigma_1$ and $\sigma_2$, respectively, intersect at exactly one point if $\sigma_1 = \sigma_2$, namely at $x = (m_1 + m_2)/2$; if $\sigma_1 \neq \sigma_2$ the density functions 
 intersect at exactly two points 
\begin{equation*}
x_{1,2} = \frac{\sigma^2_1 m_2 - \sigma^2_2 m_1 \pm \sigma_1\sigma_2 \sqrt{(m_1 - m_2)^2 + (\sigma^2_1 -\sigma^2_2) (\log \sigma^2_1 - \log\sigma^2_2)}}{\sigma^2_1 - \sigma^2_2}.
\end{equation*}
\end{proposition}

\noindent
N.B. For brevity, the standard notation $N(m, \sigma^2)$ will be used throughout this note to denote a normal distribution with mean $m$ and standard deviation $\sigma >0$. 

\begin{example}
\label{ex7feb-1}
\begin{enumerate}[label=(\roman*)]
\item
Let $P_1 \sim N(100, 10^2)$ and $P_2 \sim N(110, 10^2)$. By Proposition~\ref{prop4feb-1}, the unique crossing point of the density functions of $P_1$ and $P_2$ is at $x = 105$, which implies that the proportion of $P_2$ that is above any cutoff $c > 105$ is greater than the proportion of $P_1$ above $c$.  Conversely, the proportion of $P_1$ below any $c < 105$ is greater than the proportion of $P_2$ below $c$.
\item
Let  $P_1 \sim N(100, 10^2)$, and  $P_2 \sim N(101,11^2)$. By Proposition~\ref{prop4feb-1}, the two crossing points of the density functions of $P_1$ and $P_2$ are at $x_1 \cong 83.52$ and $x_2 \cong 106.95$, which implies that the proportion of $P_2$ that is above any cutoff $c>x_2$ is greater than the proportion of $P_1$ above $c$.  Similarly, in this case $P_2$ also dominates $P_1$ in the lower tail in that the proportion of $P_2$ that is below any cutoff $c < x_1$ is also greater than the proportion of $P_1$ below $c$.
\end{enumerate}
\end{example}

\section{Right-Tail Dominance}

As follows immediately from Proposition~\ref{prop4feb-1} and was illustrated in Example~\ref{ex7feb-1}, given any two different normal distributions, one of them necessarily dominates the other in the right tail. A much stronger conclusion is true, and that is the purpose of this section. 

Recall that a probability measure $P$ on the real line is uniquely determined by its \textit{complementary cumulative distribution function} (ccdf) $\bar{F_P}$, defined by $\bar{F_P}(x) = P((x, \infty))$ for all $x \in \mathbb{R}$.
($\bar{F}_{P}$ is also often called the \textit{survival function} of $P$, since $\bar{F}_{P}(c)$ represents the $P$-probability of the set above the cutoff threshold $c$, i.e., the fraction that survives when all values less than or equal to $c$ are removed.)

A distribution $P_1$ may be said to dominate another distribution $P_2$ in the right tail if the proportion of $P_1$ that is above a cutoff $c$ is strictly larger than the proportion of $P_2$ above $c$ for all sufficiently large $c$, i.e., in terms of the ccdf's, if
$\bar {F}_{P_1} (c) > \bar {F}_{P_2} (c)$ for all sufficiently large $c$.  A much stronger notion of domination in the right tail is if the survival probabilities of $P_1$ eventually become arbitrarily larger than those of $P_2$ as the cutoff increases;  this is formalized in the next definition.

\begin{definition} \label{dfn1}
A probability distribution $P_1$ \textit{strongly dominates distribution $P_2$ in the right tail} if 
\begin{equation*}
\lim_{c \to \infty} \frac{\bar{F}_{P_2} (c)} {\bar{F}_{P_1} (c) } = 0, 
\end{equation*}

\vspace{0.5em}
\noindent
that is, if the 
relative proportion of $P_1$ in the right tail compared to that of $P_2$ in the right tail approaches 100\% as the cutoff threshold $c$ gets arbitrarily large.
\end{definition}

Recall that a continuous (absolutely continuous) probability distribution $P$ is uniquely determined by its probability density function $f_P : \mathbb{R} \to [0,\infty)$ via $P((a,b)) = \int_a ^b f_P (x) dx$, so in particular, $\bar{F_P}(c) = \int_c^\infty f_P (x) dx$.
 The next lemma records a simple relationship between the quotients of probability density functions (pdf's) and the quotients of the corresponding ccdf's, and will be used in several examples and proofs below.

\begin{lemma}
\label{lemma1}
Suppose $P_1$ and $P_2$ are continuous probability distributions with strictly positive pdf's $f_1$ and $f_2$ and with ccdf's $\bar{F}_1$ and $\bar{F}_2$, respectively.  If 
$
\lim_{x \to \infty} {f_1 (x)}/{f_2 (x)} = \alpha$, 
then
$ \lim_{c \to \infty} {\bar{F}_1(c)}/{\bar{F}_2(c))} = \alpha$.
\end{lemma}

\begin{proof}
Straightforward from the definition of density functions and from the linearity and order-preserving properties of integration.
\end{proof}

\vspace{0.5em}
\begin{example} \label{ex1}
\begin{enumerate}[label=(\roman*)]
\item 
Let $P_1$ and $P_2$ be Cauchy distributions with medians $m_1 = 0$ and $m_2 = 0.5$ and scale parameters $s_1 = 1$ and $s_2 = 0.5$, respectively.  Then by Lemma~\ref{lemma1}, 
\begin{equation*}
\lim_{x \to \infty} \frac{f_{P_1}(x)}{f_{P_2}(x)} = 
\lim_{x \to \infty} \frac{\pi(1+2x^2 - 2x)}{\pi (1+x^2)} = 2 =
\lim_{c \to \infty} \frac{\overline{F}_{P_1}(c)}{\overline{F}_{P_2}(c)},
\end{equation*}
which implies that, as $c \to \infty$, the $P_1$-probability of the set of numbers greater than $c$ approaches exactly twice the $P_2$-probability of numbers greater than $c$.  Thus $P_1$ dominates $P_2$ in the right tail, but does not \textit{strongly} dominate $P_2$ in the right tail.

\item 
Let $P_1$ and $P_2$ be Laplace distributions with medians $m_1 = 1$ and $m_2 = 0$ and scale parameters $s_1 = s_2 = 1$, respectively.   Then 
\begin{equation*}
\lim_{c \to \infty} \frac{\bar{F}_{P_1} (c)} { \bar{F}_{P_2} (c)} = \lim_{c \to \infty} \frac{e^{1-c}}{e^{-c}}  = e,
\end{equation*}
so neither $P_1$ nor $P_2$ strongly dominates the other in the right tail.

\item 
Let $P_1$ and $P_2$ be Laplace distributions with medians $m_1 = m_2 = 0$ and scale parameters $s_1 = 1$ and $s_2 = 0.5$, respectively.   Then $P_1$ strongly dominates $P_2$ in the right tail since
\begin{equation*}
\lim_{c \to \infty} \frac{\bar{F}_{P_2} (c)} { \bar{F}_{P_1} (c)} = \lim_{c \to \infty} \frac{e^{-2c}}{e^{-c}}  = 0.
\end{equation*}

\item 
Let $P_1$ and $P_2$ be normal distributions with identical variances $+1$, and with means 1 and 0, respectively. 
Then the density functions $f_{P_1} (x) = (1/\sqrt{2\pi}) e^{-{(1/2)}(x-1)^2}$ and $f_{P_2} (x) = (1/\sqrt{2\pi}) e^{-(1/2)x^2}$ for $P_1$ and $P_2$, respectively, satisfy
$(f_{P_1}(x))/(f_{P_2}(x)) = e^{x - (1/2)} \to \infty$ as $x \to \infty$, so by Lemma~\ref{lemma1},  $P_1$ strongly dominates $P_2$ in the right tail.

\end{enumerate}
\end{example}

As was seen in Example~\ref{ex1}, for two given different Cauchy distributions or two different Laplace distributions, neither distribution may strongly dominate the other in the right tail. This is in sharp contrast to the main conclusion in this note, where it will be shown that in every finite mixture of different \textit{normal} distributions, there is always a unique one of those distributions that strongly dominates every one of the other distributions in the right tail.
The next theorem, the key step in this conclusion, follows easily from elementary basics of normal distributions; since the authors know of no explicit reference to these conclusions, a proof is included. \\

\begin{theorem} \label{prop1}
Let $P_1$ and $P_2$ be different normal distributions.  Then
\begin{enumerate}[label=(\roman*)]
\item either $P_1$ strongly dominates $P_2$ in the right tail or $P_2$ strongly dominates $P_1$ in the right tail;
\item if $P_1$ strongly dominates $P_2$ in the right tail, then either $P_1$ has greater mean (average value) than $P_2$ or $P_1$ has greater variance than $P_2$, or both;
\item if $P_1$ has greater variance than $P_2$, then $P_1$ strongly dominates $P_2$ in both right and left tails, independent of the means.
\end{enumerate}
\end{theorem}

\begin{proof}

Suppose that $P_1 \sim N(m, \sigma_1^2)$ and $P_2 \sim N(m, \sigma_2^2)$ are normal distributions with pdf's $f_1$ and $f_2$, respectively.  Since $P_1$ and $P_2$ are different, either ${\sigma_1} \neq {\sigma_2}$, or ${\sigma_1} = {\sigma_2}$ and $m_1 \neq m_2$.

\noindent
\textit{Case 1}.  $\sigma_1 \neq \sigma_2$.  Without loss of generality, $\sigma_1 > \sigma_2 = 1$.  Then
\begin{align}\label{eq1}
\frac{f{_1}(x)}{f{_2}(x)} & = \frac{\sqrt{2 \pi}}{\sigma_1 \sqrt{2 \pi}} \exp \left(  \frac{-(x-m_1)^2}{2 {\sigma_1}^2} + \frac{(x-m_2)^2}{2} \right) \\ 
& = \frac{1}{\sigma_1} \exp \left( \frac{x^2}{2}\left(1 - \frac{1}{{\sigma_1}^2}\right) + x \left(\frac{m_1}{{\sigma_1}^2} - m_2 \right) + \left(  \frac{{m_1}^2}{2} -   \frac{{m_1}^2}{2{\sigma_1}^2}  \right) \right).
\nonumber
\end{align}

Since $\sigma_1^2 > 1$, $(1 - (1/\sigma_1^2)) >0$, so (\ref{eq1}) implies that $\lim_{x \to \infty} (f{_1}(x)/f{_2}(x)) = \infty$, and thus by Lemma~\ref{lemma1}, $P_1$ strongly dominates $P_2$ in the right tail.

\vspace{1em}
\noindent
\textit{Case 2}. ${\sigma_1}^2 = {\sigma_2}^2$ and $m_1 \neq m_2$.  Without loss of generality, $\sigma_1 = \sigma_2 = 1$. Then

\begin{align}\label{eq2}
\frac{f{_1}(x)}{f{_2}(x)} & = \exp \left( \frac{-x^2 + 2m{_1}x - {m_1}^2 + x^2 - 2m{_2}x + {m_2}^2}{2} \right) \\
& = \exp \left( (m_1 - m_2) x + \frac{({m_2}^2 - {m_1}^2)}{2} \right). \nonumber
\end{align}

Also without loss of generality, $m_1 > m_2$, in which case (\ref{eq2}) implies, via Lemma~\ref{lemma1} as in Case 1, that $P_1$ strongly dominates $P_2$ in the right tail. 
This concludes the proof of (i); the proofs of (ii) and (iii) follow similarly.
\end{proof}

The same essential argument extends easily to show that among every finite collection of different normal distributions, strong domination in the right tail by one of those distributions is inevitable.

\begin{corollary}
\label{cor1}
Given a finite number of different normal distributions $P_1, \ldots, P_n$, there is a unique one of these distributions that strongly dominates all the others in the right tail.
\end{corollary}

\begin{proof}
For each $i \in \{1, \ldots, n\}$, let normal distributions $P_i$ have mean $m_i$ and standard deviation $\sigma _i$.  Since the distributions are all different, if 
$m_i = m_j$ and $\sigma_i = \sigma_j$ then $i=j$,
which implies that there exists a unique $i^* \in \{1, \ldots,n\}$ such that $m_{i^*} = \max\{m_j : \sigma_j = \max \{ \sigma_1, \ldots, \sigma_n\}\}$. By the arguments for Cases 1 and 2 in Theorem~\ref{prop1}, $P_{i^*}$ strongly dominates $P_i$ in the right tail for all $i \neq i^*$.
\end{proof}

 \begin{figure}[!ht] 
  \center\includegraphics[width=0.9\textwidth]{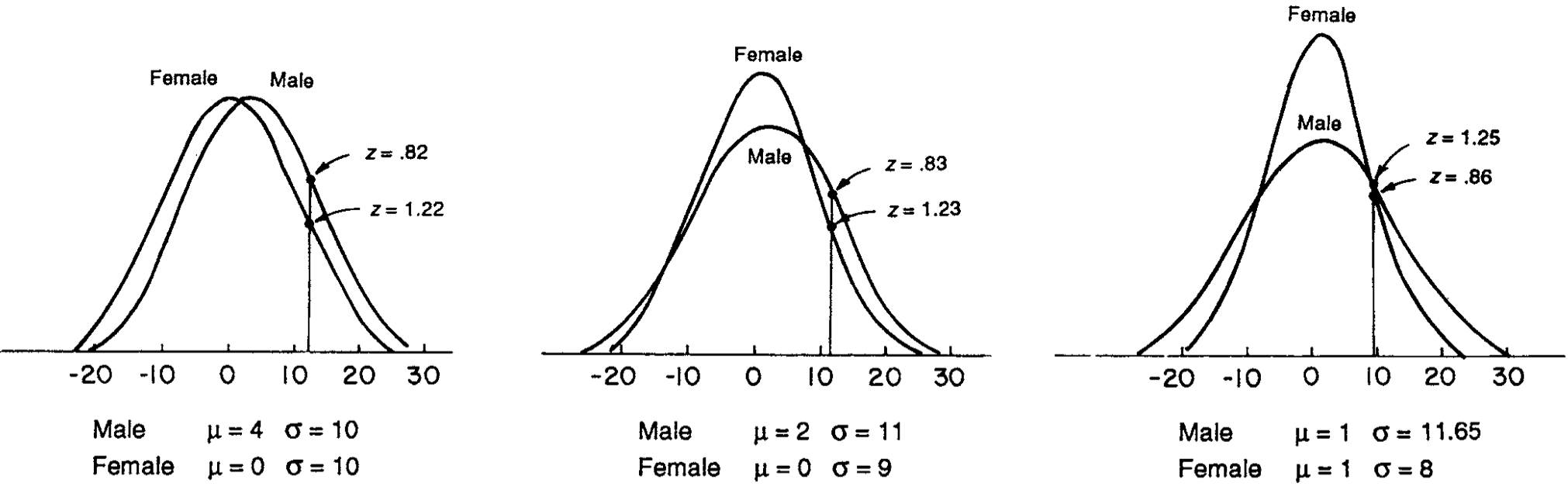}
  \caption{Feingold's ``Hypothetical superimposed male and female distributional curves for three comparisons between the sexes'' \citep[p.~7]{vob59}. Note the strong dominance in the right tail of the M subpopulation. Note also that the important case where M has smaller mean and larger variance is missing, and this is illustrated in Figure~\ref{meanVar1}(left).}
  \label{feingold2}
\end{figure} 

As was seen in Theorem~\ref{prop1}, if $P_1$ has either greater variance than $P_2$, or the same variance and higher mean, then $P_1$ will strongly dominate $P_2$ in the right tail ; see Figure~\ref{feingold2}.  Moreover, as
pointed out by Feingold, ``Most important, what might appear to be trivial group differences in both variability and central tendency can cumulate to yield very appreciable differences between the groups in numbers of extreme scorers" \cite[p.~11]{vob59}. 
The next example, a slight modification of the numerical example suggested by Feingold,  illustrates this observation with two normal distributions that are close in mean value (100 vs. 101) and in standard deviation (10 vs. 11).

\begin{figure}[!ht] 
  \center\includegraphics[width=0.9\textwidth]{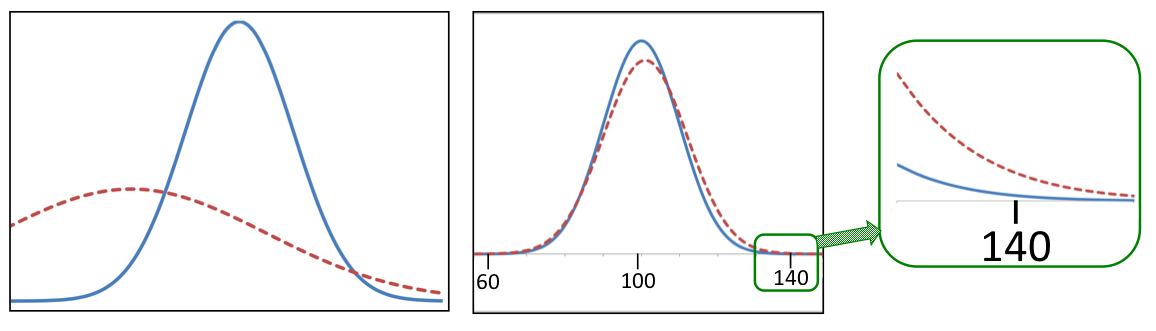}
  \caption{Comparisons of the density functions of two different normal distributions. The graph on the left complements Feingold's illustration in Figure~\ref{feingold2} above, and illustrates how even a distribution with smaller mean will eventually dominate in the right tail if its variance is larger.  The graph on the right illustrates the situation in Examples~\ref{ex7feb-1} and \ref{ex2}, with the blue (solid) curve having mean 100 and standard deviation 10, and the red (dotted) curve with mean 101 and standard deviation 11, showing how even small increases in means and standard deviations can result in significant differences in the (right) tail. }
  \label{meanVar1}
\end{figure}

\begin{example} \label{ex2}
Suppose a population $X$ consists of two mutually exclusive subpopulations $X_1$ and $X_2$, where the values of a given trait are normally distributed with distributions $P_1 \sim N(100, 10^2)$ and  $P_2 \sim N(101,11^2)$, respectively, as shown in Figure~\ref{meanVar1}(center).  Assuming an equal number of each group in the overall population $X$, routine calculations with standard \textit{erf} function programs and Newton's method yield that among the top one tenth of 1\% of the combined population, approximately 82\% of the individuals  will be type $X_2$. (To put this in perspective, 0.1\% of the current population of India is more than 1.3 million people.) 
\end{example}
\noindent
Note: The results in Theorem~\ref{prop1} and Corollary~\ref{cor1} for normal distributions only depend on first-order asymptotic terms, and the question of which more general classes of distributions with rapidly decaying tails satisfy the same conclusions is left to the interested reader.

\section{Over-Representation in the Right Tail}
\label{sec:overRep}
Whether a particular subpopulation is \textit{over-represented} or \textit{under-represented} among the other subpopulations with respect to given values for a given trait depends on the relative size of that subpopulation with those trait values compared to the size of the whole population with those trait values. For example, if subpopulation $X_1$ comprises 30\% of the total population, but comprises 40\% of the population with trait values above a given cutoff $c$, then $X_1$ is over-represented in the portion of the total population with values greater than $c$. 

The goal of this section is to show that a simple consequence of Corollary~\ref{cor1} is that in every finite mixture of different normal distributions, exactly one of those distributions will be strongly over-represented in the right tail. (Recall that a \textit{finite mixture of distributions} is a probability distribution  with cdf $F$ satisfying $F = \sum\limits_{i=1}^n w_i F_i$, where $n>1$, $F_1, \ldots, F_n$ are cdfs, and $w_1, \ldots, w_n$ are strictly positive weights with $\sum\limits_{i=1}^n w_i =1$. )

\begin{definition}
\label{dfn2}

Given a finite mixture of distributions $F = \sum\limits_{i=1}^n w_i F_i$, distribution $F_{i^*}$ is strongly over-represented in the right tail of $F$ if, as $c \to \infty$,
  the proportion of subpopulation $F_{i^*}$  with values above $c$ approaches 100\% of the total population of $F$ with values above $c$ , that is, if
\begin{equation*}
\lim_{c \to \infty} \frac{w_{i^*} \bar{F}_{i^*}(c)}{\sum\limits_{i=1}^n w_i \bar{F}_i (c)} = 1.  
\end{equation*}
\end{definition}

\begin{theorem}
\label{prop4feb-2}
In every finite mixture of different normal distributions $F = \sum\limits_{i=1}^n w_i F_i$, there is a a unique $i^* \in \{1, \ldots, n\}$ such that $F_{i^*}$ is strongly overrepresented in the right tail of $F$.
\end{theorem}

\begin{proof}
Immediate from Corollary~\ref{cor1} and the definitions of strongly dominated and strongly overrepresented.
\end{proof}

\begin{example}
\label{ex10feb-1}
The College Board SAT scores of males and females are usually assumed (or designed) to be approximately normally distributed \citep{vob198}. Unless the distributions are identical, Theorem~\ref{prop4feb-2} implies that exactly one of those two genders must be strongly over-represented in the right tail, and that this over-representation will increase as the score range increases; Figure~\ref{sawyer1} illustrates this with actual College Board SAT statistics. 
 \begin{figure}[!ht] 
 \center\includegraphics[width=0.9\textwidth]{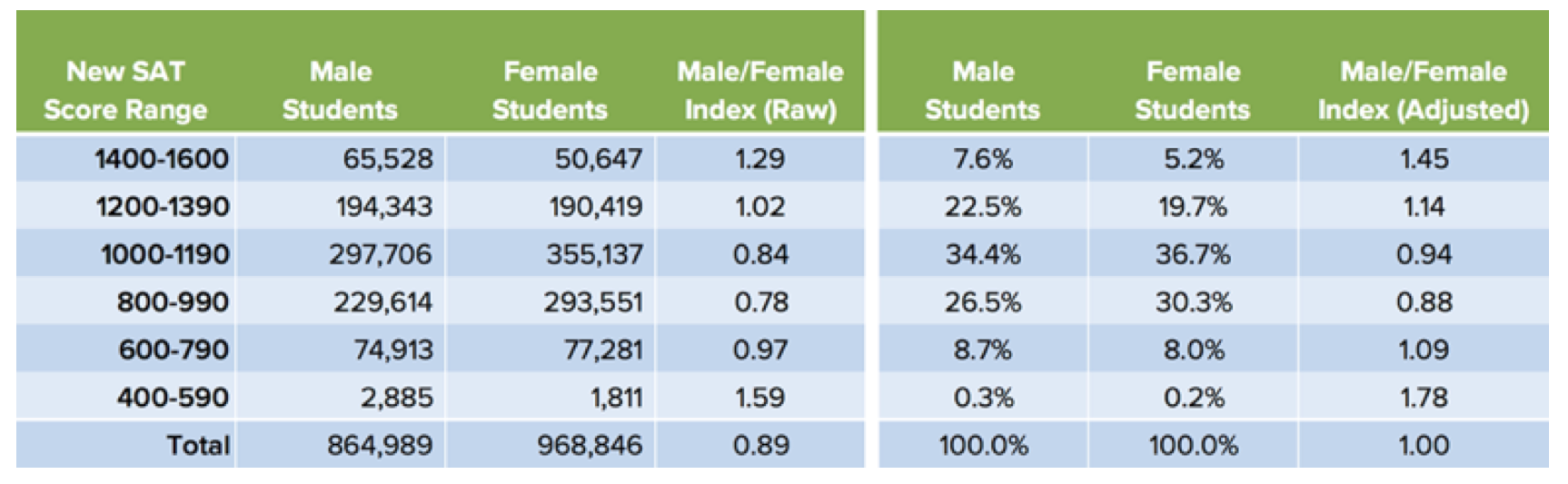}
  \caption{Statistics for nearly two million students for the 2016 Edition of the Scholastic Aptitude Test (SAT), with breakdown by gender and score ranges \citep{vob186}. Note that the proportions of males in various score ranges increase as the score range increases and also increase as the score range decreases.
About 10\% more females participated than males, which is reflected in the Adjusted Male/Female Ratios.}
  \label{sawyer1}
\end{figure} 
\end{example}

In real life examples, of course, there are no variables that are exactly normally distributed, since the normal distribution is continuous, and the numbers of people in various categories, for example, is necessarily finite. But if distributions are close to being normally distributed, the right-tail over-representation of a unique subpopulation predicted by Theorem~\ref{prop4feb-2} (and the analogous conclusions for left tails) are perhaps reasonable to expect.  Similarly, in real-life examples, calculations involving tails that are 6 or 7 standard deviations out involve probabilities of less than one in 10 billion, and are meaningless among the current human population of this planet. \\



\begin{thebibliography}{9}
\newcommand{\enquote}[1]{``#1''}
\providecommand{\natexlab}[1]{#1}
\providecommand{\url}[1]{\texttt{#1}}
\providecommand{\urlprefix}{URL }

\bibitem[{Anand and Winters(2018)}]{novob18}
Anand, Rohini and Winters, Mary-Frances.
\newblock \enquote{A Retrospective View of Corporate Diversity Training From
  1964 to the Present.}
\newblock \emph{Academy of Management Learning \& Education} 7.3 (2018):
  356--372.
\newblock
  \href{http://dx.doi.org/10.5465/AMLE.2008.34251673}{Doi:10.5465/AMLE.2008.34251673}.

\bibitem[{Ceci et~al.(2009)Ceci, Williams, and Barnett}]{vob19}
Ceci, Stephen~J., Williams, Wendy~M., and Barnett, Susan~M.
\newblock \enquote{Women's underrepresentation in science: Sociocultural and
  biological considerations.}
\newblock \emph{Psychological Bulletin} 135.2 (2009): 218--261.
\newblock \href{http://dx.doi.org/10.1037/a0014412}{Doi:10.1037/a0014412}.

\bibitem[{Dorans(2002)}]{vob198}
Dorans, Neil~J.
\newblock \enquote{Recentering and Realigning the {SAT} Score Distributions:
  How and Why.}
\newblock \emph{Journal of Educational Measurement} 39.1 (2002): 59--84.


\bibitem[{Feingold(1995)}]{vob59}
Feingold, Alan.
\newblock \enquote{The additive effects of differences in central tendency and
  variability are important in comparisons between groups.}
\newblock \emph{American Psychologist} 50.1 (1995): 5--13.
\newblock
  \href{http://dx.doi.org/10.1037/0003-066X.50.1.5}{Doi:10.1037/0003-066X.50.1.5}.

\bibitem[{Krebs et~al.(2020)Krebs, Narahari, Cook-Armstrong, Chandrabhatla,
  Mehaffey, Upchurch, and Showalter}]{novob34}
Krebs, Elizabeth~D, Narahari, Adishesh~K, Cook-Armstrong, Ian~O, Chandrabhatla,
  Anirudha~S, Mehaffey, J~Hunter, Upchurch, Gilbert~R, and Showalter, Shayna~L.
\newblock \enquote{The Changing Face of Academic Surgery: Overrepresentation of
  Women among Surgeon-Scientists with {R01} Funding.}
\newblock \emph{Journal of the American College of Surgeons.} 231.4 (2020):
  427--433.
\newblock
  \href{http://dx.doi.org/10.1016/j.jamcollsurg.2020.06.013}{Doi:10.1016/j.jamcollsurg.2020.06.013}.

\bibitem[{Lett et~al.(2019)Lett, Murdock, Orji, Aysola, and Sebro}]{novob35}
Lett, Elle, Murdock, H.~Moses, Orji, Whitney~U., Aysola, Jaya, and Sebro,
  Ronnie.
\newblock \enquote{Trends in Racial/Ethnic Representation Among {US} Medical
  Students.}
\newblock \emph{JAMA Network Open} 2.9 (2019).
\newblock
  \href{http://dx.doi.org/10.1001/jamanetworkopen.2019.10490}{Doi:10.1001/jamanetworkopen.2019.10490}.

\bibitem[{Marum(2020)}]{novob36}
Marum, Rob~J.
\newblock \enquote{Underrepresentation of the elderly in clinical trials, time
  for action.}
\newblock \emph{British journal of clinical pharmacology : BJCP.} 86.10 (2020):
  2014--2016.
\newblock \href{http://dx.doi.org/10.1111/bcp.14539}{Doi:10.1111/bcp.14539}.

\bibitem[{Roser et~al.(2019)Roser, Appel, and Ritchie}]{vob191}
Roser, Max, Appel, Cameron, and Ritchie, Hannah.
\newblock \enquote{Human Height.}
\newblock \url{https://ourworldindata.org/human-height}, 2019.
\newblock Last accessed January 22, 2021.

\bibitem[{Sawyer(2017)}]{vob186}
Sawyer, Art.
\newblock \enquote{How the New {SAT} has Disadvantaged Female Testers.}
\newblock \url
  {https://www.compassprep.com/new-sat-has-disadvantaged-female-testers/},
  2017.
\newblock Accessed January 15, 2021.

\end{thebibliography}
\end{document}